\documentclass[12 pt]{article}
\usepackage[utf8]{inputenc}

\title{Block fusion systems over maximal nilpotency class groups}
\author{Afaf Jaber}
\date{}

\usepackage[T1]{fontenc}
\usepackage{amssymb,amstext,amsmath,amsthm}
\usepackage{graphicx}
\graphicspath{{Figures/}}
\usepackage[caption=false]{subfig}
\usepackage{empheq}
\usepackage{multirow}

\usepackage{hyperref}
\usepackage{threeparttable}
\usepackage{tabularx}
\usepackage{framed}
\usepackage{algorithmic}
\usepackage{lineno}
\usepackage{tikz}
\usepackage{tikz-cd}
\usepackage[linesnumbered,ruled,vlined,longend]{algorithm2e}
\usetikzlibrary{shapes.geometric,arrows}
\usepackage{verbatim}
\usepackage{mathtools}
\newcommand\myeq{\stackrel{\mathclap{\normalfont\mbox{def}}}{=}}
\usepackage{extarrows}

\usepackage{latexsym}
\usepackage{babel}

\usepackage{background}
\backgroundsetup{
  position=current page.west,
  angle=90,
  nodeanchor=west,
  vshift=-5mm,
  opacity=1,
  scale=3,
  contents=Draft
}

\usetikzlibrary{positioning,shadows,arrows}

\usepackage{booktabs}
\usepackage{float}
\restylefloat{table}

\usepackage{lipsum}

\theoremstyle{plain}
\newtheorem{theo}{Theorem}[section]
\newtheorem{lemm}[theo]{Lemma}
\newtheorem{prop}[theo]{Proposition}
\theoremstyle{definition}
\newtheorem{defi}{Definition}[section]

\theoremstyle{remark}
\newtheorem{rema}{Remark}

\newcommand{\Hom}{{\rm Hom}}
\newcommand{\Inj}{{\rm Inj}}
\newcommand{\Aut}{{\rm Aut}}
\newcommand{\Inn}{{\rm Inn}}
\newcommand{\Out}{{\rm Out}}
\newcommand{\Ff}{\mathcal{F}}

\newcommand{\C}{{\rm C}}
\newcommand{\SL}{{\rm SL}}
\newcommand{\GL}{{\rm GL}}
\newcommand{\PGL}{{\rm PGL}}

\DeclareUnicodeCharacter{0301}{\'{e}}

\newenvironment{tightcenter}{%
  \setlength\topsep{0pt}
  \setlength\parskip{0pt}
  \begin{center}
}{%
  \end{center}
}

\title{Block fusion systems over maximal nilpotency class $3$-groups}
\author{Afaf Jaber}

\begin{document}
\maketitle
\begin{abstract}
    We generalize the Reduction Theorem of Kessar-Stancu \cite{KS} so it can be applicable to exotic fusion systems over the maximal nilpotency class of rank two $3$-groups. This is an essential step towards proving that these fusion systems are also block exotic. 
\end{abstract}

\section{Introduction}
 Let $p$ be a prime. A fusion system over a $p$-group $P$ is a category with objects the subgroups of~$P$ and morphisms satisfying certain axioms that are mentioned in Definition~$\ref{Definition of fusion systems}$.

Classical examples of fusion systems for a prime $p$ are the fusion systems of finite groups and fusion the systems of blocks of a finite group algebra over a field of characteristic $p$.

Brauer third main theorem gives the inclusion of the fusion systems of groups into the fusion systems of blocks, stating that any fusion system of a finite group is the fusion system of a the principal block of the same finite group. The opposite inclusion is an open problem. An equivalent reformulation of this problem, is showing that exotic fusion systems (i.e. general fusion systems that do not arise as fusion systems of a finite group $G$) don't appear as fusion systems of blocks of finite groups. 

This has been proven to be true by Kessar \cite{K} for Solomon exotic fusion systems at $p=2$, by Kessar-Stancu \cite{KS} for Ruiz-Viruel exotic fusion system at $p=7$ and by Serwene \cite{S} for Park-Semeraro exotic fusion systems at $p=7$.

In this article we are dealing with some families of exotic fusion systems over maximal nilpotency class $3$-groups of rank two. To our knowledge, this is the first case that was treated in the literature at $p=3$.

One cannot apply the Reduction Theorem of Kessar-Stancu \cite{KS} on these groups as the fusion systems in the families considered in this article are not reduction simple. To circumvent we generalize the reduction theorem of Kessar-Stancu by : 
\begin{theo}\label{Alternative reduction theorem}
Let $\mathcal{F}_{1}$ and $\mathcal{F}_{2}$ be fusion systems over a $p$-group $P$ such that $\mathcal{F}_{1} \subseteq \mathcal{F}_{2}$.Assume that
\begin{enumerate}
\item if $\mathcal{F}$ is a fusion system over $P$ containing $\mathcal{F}_{1}$, then $\mathcal{F} = \mathcal{F}_{1}$ or $\mathcal{F} = \mathcal{F}_{2}$,
\item if $\mathcal{F}$ is a non-trivial normal subsystem of $\mathcal{F}_{2}$ over $P$, then $\mathcal{F} = \mathcal{F}_{1}$ or $\mathcal{F} = \mathcal{F}_{2}$.
\item  if $\mathcal{F}$ is a $\mathcal{F}_{1}$ or $\mathcal{F}_{2}$ fusion system, $N$ is a proper non trivial strongly $\mathcal{F}$-closed subgroup of $\mathcal{F}$, and $\mathcal{\tilde{F}}$ is a fusion system containing $\mathcal{F}$, then there is no normal fusion subsystem of $\mathcal{\tilde{F}}$ over $N$.
\end{enumerate}
If there exists a finite group with an $\mathcal{F}_{1}$ or $\mathcal{F}_{2}$-block, then there also exists a quasisimple group $L$ with $p'$-centre having an $\mathcal{F}_{1}$ or $\mathcal{F}_{2}$-block.
\end{theo}

We conclude by showing that the considered families of fusion system over the maximal nilpotency class $3$-groups,  $B(3,4;0,2,0)$ and $B(3,2k,0,\gamma,0)$ for $k>2$ and $\gamma\ne 0$, satisfy the conditions of the Theorem \ref{Alternative reduction theorem} using the following two theorems. The definitions of the $3$-groups $B(3,4;0,2,0)$ and $B(3,2k,0,\gamma,0)$ for $k>2$ and $\gamma\ne 0$ will be given in Section 4.

\begin{theo}\label{Condition for ART B(3,4,0,2,0)}
Let $B:=B(3,4;0,2,0)$ and $\mathcal{F}$ be an exotic fusion system over $B$.
Let $\mathcal{\tilde{F}}$ be a fusion system containing $\mathcal{F}$, $T$ be an $\mathcal{F}$-strongly closed subgroup of $B$ and $\mathcal{F'}$ be a fusion system of $\mathcal{\tilde{F}}$, over $T$.
Then $\Aut_{\mathcal{F'}}(T)$ is an index $3$ subgroup of $\Aut_{\tilde \Ff}(T)$ and $\Aut_{\mathcal{F'}}(T)\ntrianglelefteq \Aut_{\tilde \Ff}(T)$. Hence $\mathcal{F'}\ntrianglelefteq\mathcal{\tilde{F}}$.
\end{theo}

\begin{theo}\label{Condition for ART B(3,2k,0,gamma,0)}
Let $B:=B(3,2k,0,\gamma,0)$ for $ \gamma\neq 0,\, k\geq 3$ and $\mathcal{F}$ be an exotic fusion system over $B$.
Let $\mathcal{\tilde{F}}$ be a fusion system containing $\mathcal{F}$, $T$ be an $\mathcal{F}$-strongly closed proper non trivial subgroup of $B$ and $\mathcal{F'}$ be the fusion system of $\mathcal{\tilde{F}}$, over $T$.
Then there exists $\gamma \in \Aut_{\mathcal{F}}(T)$ such that $^{\gamma}\Aut_{\mathcal{F'}}(T)\neq \Aut_{\mathcal{F'}}(T)$. Hence $\mathcal{F'}\ntrianglelefteq\mathcal{\tilde{F}}$
\end{theo}

Hence, if there is an exotic fusion system over these two groups that is block realisable, then there exists a quasisimple group $L$ with the same fusion system structure. \\

To show that the considered exotic fusion system are also block-exotic, one needs to show that there is no quasisimple group having a block with that fusion structure.
This is still work in progress.

\section{Fusion Systems and basic properties}
We start by giving the definition of fusion systems. For the classical results on fusion systems we refer the reader to \cite{AKO}.

\begin{defi} \label{Definition of fusion systems}
A fusion system over a finite $p$-group $P$ is a category $\mathcal{F}$, with objects the
set of all subgroups of $P$, and,  all $R$ and $Q$ subgroups of $P$, sets of morphisms Hom$_{\mathcal{F}}(R,Q)$ that satisfy the following two properties :
\begin{enumerate}
\item Hom$_{P}(Q, R) \subseteq$ \Hom$_{\mathcal{F}} (Q, R) \subseteq$ \Inj$(Q, R)$. \\where Hom$_{P}(Q, R)$ are just the conjugation maps from $Q$ in $R$ by elements of $P$ and Inj$(Q, R)$ are the injective maps from $Q$ in $R$.
\item each $\phi \in$ Hom$_{\mathcal{F}}(Q, R)$ is the composite of an isomorphism in $\mathcal{F}$ followed by inclusion.
\end{enumerate}
\end{defi}

The classical example of a fusion system is defined by the conjugation action of a finite group $G$ over a $p$-subgroup.

\begin{defi}
A fusion system $\mathcal{F}$ of a finite group $G$ over a $p$-subgroup $P$ of $G$ is the category where:
\\-The objects are subgroups of $P$
\\-The morphisms are the conjugation maps $c_{g}$ by elements $g \in G$ between subgroups $Q$ and $R$ of $P$ :  Hom$_{\mathcal{F}}(Q,R)=\{c_{g} \mid\,^{g}Q \leq R,\,g\in G\}$.\\
This category is noted $\mathcal{F}=\mathcal{F}_{P}(G)$
\end{defi}

The above definition of fusion systems, is broad and does not capture the specific properties - given by Sylow theorems - of the case when $P$ is a Sylow $p$-subgroup of $G$. That's why we use the notion of "saturated" fusion systems.
To define the saturation property, we need to introduce some more notions.

\begin{defi}
Let $\mathcal{F}$ be a fusion system over a $p$-group $P$ and let $Q$, $R$ be subgroups of $P$.
\begin{enumerate}
\item A morphism in $\mathcal{F}$ from $P$ to $Q$ is called an $\mathcal{F}$-morphism.
\\ An isomorphism in $\mathcal{F}$ from  $Q$ to $R$ is called an $\mathcal{F}$-isomorphism. 
\item $Q$ and $R$ are called $\mathcal{F}$-conjugate if
there exists an $\mathcal{F}$-isomorphism from $Q$ to $R$. The $\mathcal{F}$-isomorphism class of $Q$ is noted $Q^{\mathcal{F}}$.\\ By analogy, the $\mathcal{F}$-isomorphism class of an element $x \in P$ is noted $x^{\mathcal{F}}=\lbrace \varphi (x) \mid \varphi \in \Hom_{\mathcal{F}}(\langle x  \rangle , P) \rbrace,$ where $\langle x \rangle$ is the subgroup of $P$ generated by $x$.
\item $Q$ is called fully
automised in $\mathcal{F}$ if $\Aut_{P}(Q)$ is a Sylow p-subgroup of $\Aut_{\mathcal{F}}(Q)$.
\item $Q \leq P$ is fully centralised in $\mathcal{F}$ if $\vert C_{P}(Q)\vert \geq \vert C_{P}(R)\vert$ for all $R \in Q^{\mathcal{F}}$ .
\item $Q \leq P$ is fully normalised in $\mathcal{F}$ if $\vert N_{P}(Q)\vert \geq \vert N_{P}(R)\vert$ for all $R \in Q^{\mathcal{F}}$ .
\item  $Q$ is called receptive in $\mathcal{F}$ if the following holds: for any $\mathcal{F}$-isomorphism $\phi : R \rightarrow Q$ there exists an $\mathcal{F}$-morphism $\overline{\phi} : N_{\phi} \rightarrow P$ whose restriction to $Q$ equals $\phi$,
where $N_{\phi} = \lbrace g \in N_{P}(R) \mid (x \rightarrow \phi(\phi^{-1}(x)^{g})) \in \Aut_{P}(Q)\rbrace.$

\end{enumerate}
\end{defi}

Now we can define the saturated fusion systems:
\begin{defi} 
A saturated fusion system $\mathcal{F}$ is a fusion system such that every $\mathcal{F}$-isomorphism class contains a subgroup that is fuly automised and receptive in $\mathcal{F}$.
\end{defi}

The classical example of a saturated fusion system over a $p$-group is the fusion system  $\mathcal{F}_{P}(G)$ for $P$ a Sylow $p$-subgroup of $G$.


\begin{defi}
We say that $\mathcal{F}$ is group realisable, when $\mathcal{F}=\mathcal{F}_{P}(G)$ Where $G$ is a finite group and $P$ is a Sylow $p$-subgroup $G$. Othewise we say that $\mathcal{F}$ is an exotic fusion system. 
\end{defi}
Before we start talking about the block fusion systems, we define some important classes of subgroups for a given fusion system.

\begin{defi}\label{a group is normal, centric radical etc}
Let $\mathcal{F}$ be a fusion system over a $p$-group $P$.
\begin{enumerate}

\item A subgroup $Q$ of $P$ is $\mathcal{F}$-centric ($Q \in \mathcal{F}^{c}$) if $C_{P}(R) = Z(R)$ for all $R \in Q^{\mathcal{F}}$. Equivalently, $Q$ is $\mathcal{F}$-centric if $Q$ is fully centralized in $\mathcal{F}$ and $C_{P}(Q) = Z(Q)$. We denote the set of $\mathcal{F}$-centric subgroups of $P$ by $\mathcal{F}^{c}$.
\item A subgroup $Q$ of $P$ is $\mathcal{F}$-radical (i.e. $Q \in \mathcal{F}^{r}$) if $Out_{\mathcal{F}}(Q)$ is $p$-reduced; i.e., if $O_{p}(\Out_{\mathcal{F}}(Q)) = 1.$ We say that $Q$ is $\mathcal{F}$-centric-radical if it is $\mathcal{F}$-centric and $\mathcal{F}$-radical. We denote the set of $\mathcal{F}$-centric-radical subgroups of $P$ by $\mathcal{F}^{cr}$.
\item An $\mathcal{F}$-Alperin subgroup of $P$ is a $\mathcal{F}$-centric-radical, fully $\mathcal{F}$-normalized subgroup of $P$.

\item  $Q$ is strongly closed in $\mathcal{F}$ if no element of $Q$ is $\mathcal{F}$-conjugate to an element of $
P\setminus Q$. More
generally, for any $R \leq P$ which contains $Q$, $Q$ is strongly closed in $R$ with respect to
$\mathcal{F}$ if no element of $Q$ is $\mathcal{F}$-conjugate to an element of $R \setminus Q$.
\end{enumerate}
Note the following implications, for any $\mathcal{F}$ and any $Q \leq P$:
\begin{tightcenter}
$Q$ strongly closed in $\mathcal{F}$ $ \Rightarrow Q \unlhd P.$
\end{tightcenter}
\end{defi}

By Alperin fusion theorem a saturated fusion system over $P$ is generated by the automorphisms of the $\mathcal{F}$-Alperin subgroups of $P$.
\begin{theo}\cite{AKO}
Let $\mathcal{F}$ be a saturated fusion system over a finite $p$-group $P$. Then
\begin{center}
$\mathcal{F} = \langle\, \Aut_{\mathcal{F}} (Q )\vert\, Q\le P ,Q\text{ is }\mathcal{F}\text{-Alperin}\,\rangle$ .
\end{center}
\end{theo}

Now we introduce block fusion systems. To do that we begin by introducing Brauer pairs and the Brauer homomorphism.

\begin{defi}
 Let $G$ be a finite group. $p$ be a prime dividing the order of the group $G$ and $k$ a field of characteristic $p$.\\
A $kG$-Brauer pair is a pair $(R, e)$ where $R$ is a $p$-subgroup of $G$ and $e$ is a primitive central idempotent of $kC_{G}(R)$.
\end{defi}

Since $C_{G}(R^{g}) = C_{G}(R)^{g}$, this gives a bijection between the block idempotents of $kC_{G}(R)$ and $kC_{G}(R^{g})$. Thus, $G$ acts on the Brauer pairs by conjugation, sending $(R,e)$ to  $(R^{g},e^{g})$.

\begin{defi}
Let $R$ be a $p$-subgroup of $G$ and $kC_{G}(R)=\lbrace x \in kC_{G}(R)\vert x^{r}=x, \forall r\in R \rbrace$. Let $br_{R} : (kG)^{R} \rightarrow kC_{G}(R)$ be the map which sends an element $ \sum_{x \in G} \alpha_{x}x$ of $kG$ to the element $\sum_{x \in C_{G}(R)} \alpha_{x}x$ of $kC_{G}(R)$.

This map $br_{R}$ is clearly $k$-linear and surjective and is called the Brauer morphism.
\end{defi}

We have a partial order relation on the $kG$-Brauer pairs.

\begin{defi}
Let $G$ be a finite group, and let $Q \unlhd R$ be $p$-subgroups of $G$. If $e$ is a block idempotent of $kC_{G}(R)$, then there is a unique $R$-stable block idempotent $f$ of $kC_{G}(Q)$ (i.e. a block of $kC_{G}(Q)$ that is fixed by $R$) such that $br_{R}(f)e = e$. We write then, that $(Q,f) \unlhd (R,e)$. We define the partial order $\le$ on the $kG$-Brauer pairs as the transitive closure of $\unlhd$.

If $f_{0}$ is any other $R$-stable block idempotent of $kC_{G}(Q)$, then $br_{R}(f_{0})e = 0$.
\end{defi}

\begin{defi}
Let $G$ be a finite group, $k$ an algebraically closed field of characteristic $p$, and $b$ a block of $kG$.
\begin{enumerate}
\item A $b$-Brauer pair is a Brauer pair $(R, e)$ such that $br_{R}(b)e = e$ i.e. $(1,b) \leq (R,e)$.
\item A defect group of $b$ is $p$-subgroup $P$ of $G$ of maximal order such, that $br_{P}(b) \neq 0$. Hence, if $(P, e)$ is $b$-Brauer pair, it is a maximal one 
\end{enumerate}
\end{defi}

We can now define the fusion system of a block of $kG$.

\begin{defi}
Let $b$ be a block idempotent of $kG$, and let $(P, e)$ be a maximal $b$-Brauer pair. The fusion system of the block algebra $kGb$, denoted $\mathcal{F} = \mathcal{F}_{(P,e)}(kG, b)$, is the category:
\newline - whose objects are the subgroups of $P$
\newline - whose morphisms are given by the following. Let $Q$ and $R$ be subgroups of $P$, and $e_{Q}$ and $e_{R}$ be the unique block idempotents such that $(Q, e_{Q}) \leq (P, e)$ and $(R, e_{R}) \leq (P, e)$. The set $\Hom_{\mathcal{F}} (Q, R)$ is the set of conjugation maps, $c_{x} : Q \rightarrow R$, for all elements $x$ of $G$ such that $^{x}(Q, e_{Q}) \leq (R, e_{R})$.
\newline 
When $\mathcal{F} = \mathcal{F}_{(P,e)}(kG, b)$ we say that $\mathcal{F}$ is block realisable.
\end{defi}

Remark that $\mathcal{F}_{(P,e)}(kG, b)$ is a saturated fusion system over $P$. 

\begin{rema}\label{remark about normal subfusion of block}
It has been proven that if $N\unlhd G$, then $\mathcal{F}_{P\cap N}(N) \unlhd \mathcal{F}_{P}(G)$ for $P$, $P\cap N$ Sylow $p$-subgroups of $G$, $N$. This  doesn’t work for fusion systems of blocks in general. If $c$ is block idempotent of $kN$ with defect $P \cap N$, covered by $b$, then, in general, $\mathcal{F}_{P\cap N}(kN, c)$ is not even a fusion subsystem of $\mathcal{F}_{P}(kG, b)$.
\end{rema}

 To transfer properties from fusion systems of blocks of finite groups to fusion systems of blocks of their quasisimple subquotients one needs to have relations between $\mathcal{F}_{P\cap N}(kN, c)$ and , where the notations are as in the remark above. $\mathcal{F}_{P}(kG, b)$, where the notations are as in the remark above. One way to approach this problem is via the notion of generalized Brauer pairs :

\begin{defi}
Let $G$ be a finite group, $N$ a normal subgroup of $G$ and $b$ a block $kN$.
A $(kN,G)$-Brauer pair is a pair $(P,e_{P})$, where $P$ is a $p$-subgroup of $G$ and $e_{P}$ is a block of $kC_{N}(P)$.
\end{defi}

\begin{defi}\cite{KS}[Definition 3.1]
A $(kN,b, G)$-Brauer pair is a pair $(P, e_{P})$ where $P$ is a $p$-subgroup of $G$ such that $br^{N}_{P}(b) \neq 0$ and $e_{P}$ is a block of $kC_{N}(P)$ such that $Br^{N}_{P}(b)e_{P} \neq 0$. 
\end{defi}

Remark that, When $G = N$, a $(kN,b,G)$-Brauer pair is a $b$-Brauer pair

\begin{prop}\cite{AKO}[Proposition 4.3.13]
Let $G$ be a finite group, $N$ a normal subgroup of $G$ and $b$ a block $kG$.
\begin{enumerate}
\item $(kN,G)$-Brauer pair $(R,e)$ is a $(kN,b,G)$-Brauer pair iff $br_{R}(b)e\neq 0$.
\item Any $(kN,G)$-Brauer pair is a $(kN,b,G)$-Brauer pair for a unique primitive idempotent of $(kN)^{G}.$
\end{enumerate}
\end{prop}

\begin{defi}\cite{KS}[Definition 3.3]
Let $N$ be a normal subgroup of $G$, $b$ a block of $kN^{G}$ and $(P,e_{P})$ a maximal $(kN,b,G)$-Brauer pair. For each $Q \leq P$, let $e_{Q}$ denote the unique block of $kC_{N}(Q)$ s.t $(Q,e_{Q}) \leq (P,e_{P})$. \\
Denote $\mathcal{F}_{(P,e_{P})}(kN,b,G)$ the category whose objects are subgroups of $P$ and which has morphism set Hom$_{\mathcal{F}_{(P,e_{P})}(kN,b,G)}(Q,R)=$ Hom$_{G}((Q,e_{Q}),(R,e_{R}))$ for all $Q\leq R$. \\
If $G=N$, then $\mathcal{F}_{(P ,e_{P})}(kN,c,G)$ is the usual fusion system of the block $c$, and we denote it by $\mathcal{F}_{(P ,e_{P})}(G, c)$.
\end{defi}

By [Theorem 3.4, \cite{KS}] we get that $\mathcal{F}_{(P,e_{P})}(kN,b,G)$ is a saturated fusion system over~$P$.\\\\
Next theorem gives a solution to the problem in Remark.$\ref{remark about normal subfusion of block}$.

\begin{theo}\cite{KS}[Theorem 3.5]\label{generalized fusion system}
Let $G$ be a finite group, $N\unlhd G$, $c$ be a block of $kN^{G}$ and  $\mathcal{F}_{(D,e_{D})}(kN,c,G)$ be a fusion system over $D$. If $b$ is a block of $kG$ covering $c$ with maximal $b$-Brauer pair $(P,e_{P})$, we have $P\leq D$, 
\\$\mathcal{F}_{(P,e_{P})}(kG,b)\leq \mathcal{F}_{(D,e_{D})}(kN,c,G)$ . Moreover $P \cap N = D \cap N$ and  
\\$\mathcal{F}_{(D\cap N,e_{D\cap N})}(kN,c)\unlhd \mathcal{F}_{(D,e_{D})}(kN,c,G).$
\end{theo}

\section{Reduction Theorem}
In this section we state a reduction theorem of Kessar-Stancu. Then we prove a generalization that is necessary when dealing with fusion systems of blocks with defect groups of maximal nilpotency class.
\begin{defi}

Let $\mathcal{F}$ be a fusion system over a $p$-group $P$. If $P$ does not have any non-trivial proper strongly $\mathcal{F}$-closed subgroups, we say that $\mathcal{F}$ is reduction simple
\end{defi}





We state now the reduction theorem of Kessar-Stancu :

\begin{theo}\cite{KS}[Theorem 4.2]
Let $\mathcal{F}_{1}$ and $\mathcal{F}_{2}$ be fusion systems over a $p$-group $P$ such that $\mathcal{F}_{1} \subseteq \mathcal{F}_{2}$. Assume that
\begin{enumerate}

\item $\mathcal{F}_{1}$ is reduction simple,
\item if $\mathcal{F}$ is a fusion system over $P$ containing $\mathcal{F}_{1}$, then $\mathcal{F} = \mathcal{F}_{1}$ or $\mathcal{F} = \mathcal{F}_{2}$,
\item if $\mathcal{F}$ is a non-trivial normal subsystem of $\mathcal{F}_{2}$, then $\mathcal{F} = \mathcal{F}_{1}$ or $\mathcal{F} = \mathcal{F}_{2}$.
\end{enumerate}
If there exists a finite group with an $\mathcal{F}_{1}$ or $\mathcal{F}_{2}$-block, then there also exists a quasisimple group $L$ with $p'$-centre with an $\mathcal{F}_{1}$ or $\mathcal{F}_{2}$-block.
\end{theo}


Theorem \ref{Alternative reduction theorem} is a generalisation of the Theorem above. We finish this section by giving the proof of this generalized reduction theorem :


\begin{proof}[Proof of Theorem \ref{Alternative reduction theorem}] 
Let $G$ be a minimal order group having an $\mathcal{F}_{1}$ or an $\mathcal{F}_{2}$-block $b$. By a standard reduction (see for example [\cite{K} , Proposition 2.11]), if $N$ is a normal subgroup of $G$ and $c$ is a block of $N$ with $bc = 0$, then $c$ is $G$-stable.

By abuse of notation, $B$ is a $(b, G)$-defect group. Let $H := \lbrace ^{g}B \mid g \in G \rbrace $ be the normal subgroup of $G$ generated by all $G$-conjugates of $B$ . Let $d$ be the unique block of $kH$ covered by $b$.
Given that $d$ is $G$-stable, $G$ acts by conjugation on $kHd$. Let $N$ be the kernel of the homomorphism $G \rightarrow \Out(kHd) = \Aut(kH d)/\Inn(kH d)$. Then by [\cite{ku2}, §6 theorem] $G/N$ is a $p'$-group.
We prove, using the minimality of $G$ and the hypothesis on $\mathcal{F}_{1}$ and $\mathcal{F}_{2}$, that $G = N$.

Let $c$ be the block of $N$ covered by $b$, i.e. $bc=0$ (in this case in fact we have $c = b$). Let
$(B , e_{B} )$ be a maximal $b$-Brauer pair and let $(S, e'_{S})$ be a maximal $(G, c)$ Brauer pair as in Theorem $\ref{generalized fusion system}$. Since $G/N$ is a $p'$-group, it follows that $S = B$ . Hence, we have $\mathcal{F}_{(B ,e_{B} )}(G, b)$ is a
subsystem of $\mathcal{F}_{(B ,e'_{B} )}(G, N , c)$ and that $\mathcal{F}(B ,e'_{B} )(N ,c)$ is a normal subsystem of $\mathcal{F}(B ,e'_{B} )(G, N , c)$.

Given that $b$ is an $\mathcal{F}_{1}$- or $\mathcal{F}_{2}$-block and that $\mathcal{F}_{1}$ and $\mathcal{F}_{2}$ are the only fusion systems over $B$ that contain $\mathcal{F}_{2}$ we obtain that $\mathcal{F}(B ,e_{B} )(G, N , c)$ is either $\mathcal{F}_{1}$ or $\mathcal{F}_{2}$. Again, since the only normal proper fusion subsystem over $B$ contained in $\mathcal{F}_{1}$ is $\mathcal{F}_{2}$ and $\mathcal{F}_{2}$ has no normal fusion subsystem it follows that $\mathcal{F}(B ,e'_{B} )(N , c)$ is either $\mathcal{F}_{1}$ or $\mathcal{F}_{2}$. By the minimality of $G$ we deduce that $G = N$.

As $b$ and $d$ have the same defect group $B$ and $G$ acts on $kHd$ by inner automorphisms, using
another result of Külshammer [\cite{ku}, Theorem 7], we have that $kGb$ and $kHd$ have isomorphic
source algebras, so $c$ is also an $\mathcal{F}_{1}$- or $\mathcal{F}_{2}$-block. Thus, once again by the minimality of $G$ we have $G = H$.

Let $M$ be a proper normal subgroup of $G$. Then we have three cases for the strongly $\mathcal{F}_{1}$ (or $\mathcal{F}_{2}$)-closed subgroup $B \cap M$ of $B$ , hence $B \cap M = 1$ or $B\cap M=N$, where $N$ proper, non-trival subgroup of $B$, or $B \cap M = B$. Suppose first that $B \cap M = B$ . Then $B$ and all its $G$-conjugates lie in $M$. Thus $G = M$, which is a contradiction. 
We look then at the case of $B \cap M=N$, proper non-trivial, but by the third condition of the theorem there is no fusion system over $N$ normal in  $\mathcal{\tilde{F}}$.
We're left with the case of $B \cap M = 1$. A variation of Fong reduction allows us to deduce that there is a central $p'$-extension $G'$ of $G/M$ having an $\mathcal{F}_{1}$- or $\mathcal{F}_{2}$-block (see for example [\cite{K}, Section 3; 3.3 and below]).

\end{proof}

\section{Maximal nilpotency class rank two \texorpdfstring{$3$-}-groups.}
We start by introducing the maximal nilpotency class $3$-groups of rank two and some of their properties. Most of the results are taken from the article of Diaz-Ruiz-Viruel \cite{DRV}. \\For $x$ and $y$ elements of a group denote $[x, y] = xyx^{-1}y^{-1}$.

\begin{theo}\cite{B}\label{introduction to maximal nilpotency class groups}
The non cyclic $3$-groups of maximal nilpotency class and of order $3^{r} > 3^{3}$ are the groups $B(3, r; \beta, \gamma, \delta)$ with $(\beta, \gamma, \delta)$ taking the values:\\\\
$\bullet$ For any $r>4$, $(\beta, \gamma, \delta) = (1, 0, \delta)$, with $\delta \in  \lbrace 0, 1, 2 \rbrace$.\\\\
$\bullet$ For even $r>3$, $(\beta, \gamma, \delta) \in \lbrace (0, \gamma, 0),(0, 0, \delta)\rbrace$, with $\gamma \in \lbrace 1, 2\rbrace$ and $\delta \in \lbrace 0, 1\rbrace$.\\\\
$\bullet$ For odd $r>4$, $(\beta, \gamma, \delta) \in \lbrace (0, 1, 0),(0, 0, \delta) \rbrace$, with $\delta \in \lbrace 0, 1 \rbrace $.\\
With these parameters, $B(3, r; \beta, \gamma, \delta)$ is the group of order $3^{r}$ defined by the set of generators $\lbrace s, s_{1}, s_{2}, \cdots , s_{r-1}  \rbrace$ and the following relations :
\begin{enumerate}
\item $s_{i} = [s,s_{i-1}]$ for $i \in \lbrace 2, 3, \cdots , r-1 \rbrace$,
\item $[s_{2},s_{1}] = s_{r-1}^{-\beta}$,
\item $[s_{i},s_{1}] = 1$ for $i \in \lbrace 3, 4, . . . , r-1 \rbrace$,
\item $s^{3} = s_{r-1}^{\delta}$, 
\item $s_{3}s^{3}_{2}s^{3}_{1} = s^{\gamma}_{r-1},$
\item $s_{i+2}s^{3}_{i+1}s^{3}_{i} = 1$ for $i \in \lbrace 2, 3, \cdots , r-1\rbrace$, and assuming $s_{r} = s_{r+1} = 1$. 
\end{enumerate}
\end{theo}

We recall some properties of subgroups of $B(3, r; \beta, \gamma, \delta)$ :
\begin{prop}\cite{DRV}[Proposition A.9]
Consider the notations in Theorem $\ref{introduction to maximal nilpotency class groups}$. Then, the following hold :
\begin{enumerate}

\item $\gamma_{i}(B(3, r; 0, \gamma, 0)) = \langle s_{i}, s_{i+1}, . . . , s_{r-1} \rangle$ are characteristic subgroups of order $3^{r-i}$ generated by $s_{i}$ and $s_{i+1}$ for $i = 1, .., r - 1$ (assuming $s_{r} = 1$)\,.
\item $\gamma_{1}(B(3, r; 0, \gamma,0))$ is abelian .
\item $Z(B(3, r; 0, \gamma, 0)) = \gamma_{r-1}(B(3, r; 0, \gamma,0)) = \langle s_{r-1}\rangle\,$.
\end{enumerate}
\end{prop}

Now we mention Lemma A.10 and Lemma 5.2 of \cite{DRV} that gives us more information on the form of the normal subgroups of $B(3, r; 0, \gamma, 0)$ and the proper $\mathcal{F}$-Alperin subgroups of $B(3, r; 0, \gamma, 0)$ :

\begin{rema}
It should be noted that in what follows we restrict to $r$ even and the parameters $(\beta,\gamma,\delta)$ to $ (0,\gamma,0)$, with $\gamma \in \lbrace 1,2 \rbrace$. The exotic fusion systems we are studying are over these $3$-groups.
\end{rema}

\begin{lemm}\label{normal subgroups of B(3,r,0,gamma,0)}
Let $T$ be a nontrivial proper normal subgroup in $B(3, r; 0, \gamma, 0)$. Then
\begin{enumerate}

\item $T$ contains $Z(B(3, r; 0, \gamma, 0))$.
\item If $T$ contains $s$ then $T\cong 3^{1+2}_{+}$ if $r = 4$ or $T \cong B(3, r - 1; 0, 0, 0)$ if $r > 4$.
\end{enumerate}

\end{lemm}
We can easily compute the orders of the generators of $ :B(3, r; 0, \gamma, 0)$ :
\begin{prop} \label{order of si}
For the groups $B(3, r; 0, \gamma, 0)$, the orders of generators $s_{i}$ are $o(s_{i})= 3^{\lfloor (r+1-i)/2 \rfloor}$.
\begin{proof}
Using the fact that $\gamma_{1}(B(3, r; 0, \gamma,0))$ is abelian and equation 6. in Theorem $\ref{introduction to maximal nilpotency class groups}$, we get for $i = r-2$,  $o(s_{r-2})=3$, and for $i = r-3$, we get  $o(s_{r-3})=9$. By induction procedure, taking care of the parity of r, we get the desired result.
\end{proof}
\end{prop}

\begin{lemm}  \label{An F alperin subgroup of B(3,r;0,gamma,0) is}
Let $\mathcal{F}$ be a saturated fusion system over $B(3, r; 0, \gamma, 0)$ and let $P$ be a proper
$\mathcal{F}$-Alperin subgroup. Then $P$ is one of the following table :
\begin{table}[H]
\begin{center}
\begin{tabular}{|c|c|c|}
	\toprule
	Isomorphism type & Subgroup (up to conjugation) & Conditions \\
	\hline\hline
	$\mathbb{Z}/3^{k}\mathbb{Z} \times \mathbb{Z}/3^{k}\mathbb{Z}$ & $\gamma_{1}= \langle s_{1}, s_{2} \rangle$ & $r=2k+1$\\
	\hline
	
	\multirow{2}{*}{$3^{1+2}_{+}$} & \multirow{2}{*}{$E_{i}\myeq\langle \zeta, \zeta ', ss_{1}^{i}\rangle$} & $\zeta=s_{2}^{3^{k-1}}$, $\zeta '=s_{1}^{3^{k-1}}$ for $r=2k+1$,\\
	& &  $\zeta=s_{1}^{3^{k-1}}$, $\zeta '=s_{2}^{-3^{k-2}}$ for $r=2k$,\\
	\cline{1-2}
	\multirow{2}{*}{$\mathbb{Z}/3\mathbb{Z} \times \mathbb{Z}/3\mathbb{Z}$} & \multirow{2}{*}{$V_{i} \myeq \langle \zeta, ss_{1}^{i}\rangle$} & $i \in \lbrace -1, 0,1 \rbrace$ if $\gamma=0$ \\
	& & and $i=0$ if $\gamma=1,2$\\
	\hline
	\end{tabular}
	\end{center}
\end{table}
Moreover, for the subgroups in the table, $\gamma_{1}$ and any $E_{i}$ are always $\mathcal{F}$-centric, while any $V_{i}$ is $\mathcal{F}$-centric only if it is not $\mathcal{F}$-conjugate to $\langle \zeta, \zeta ' \rangle \cong \mathbb{Z}/3 \times \mathbb{Z}/3$
\end{lemm}

The following classification of fusion systems over $B(3, 4; 0, 2, 0)$ and \allowbreak
$B(3, 2k; 0, \gamma, 0)$ is given in Theorem 5.9 in \cite{DRV}. When there is no group in the last column, the fusion system is exotic. We address the question whether these exotic fusion system are block-exotic. 
\begin{theo}

Let $(B, \mathcal{F})$ be a saturated fusion system with at least one proper $\mathcal{F}$-Alperin subgroup with $B$ a rank two $3$-group of maximal nilpotency class for $r\geq 4$. Then it must correspond to one of the cases listed in the following tables.

$\bullet$ If $B \cong B(3, 4; 0, 2, 0)$ then the outer automorphism group of the $\mathcal{F}$-Alperin subgroups
are in the following table:

\begin{table}[H]
  \begin{center}
    \caption{Fusion Systems over $B(3,4;0,2,0)$.}
    \label{tab:table1}
    \begin{tabular}{|c|c|c|c|} 
      \toprule
      $B$ & $E_{0}$ & $V_{0}$ & group realization  \\\hline

      $<\omega>$ & - & $\SL_{2}(3)$ & $\mathcal{F}(3^{4},3)$ \\\hline
      
      $<\eta\omega>$ & $\SL_{2}(3)$ & - & $E_{0}._{2}\SL_{2}(3)$ \\\hline 
      
      \multirow{2}{*}{$<\eta,\omega>$} & - & $\GL_{2}(3)$ & $\mathcal{F}(3^{4},3).2$\\\cline{2-4}
      
      & $\GL_{2}(3)$ & - &  $E_{0}._{2}\GL_{2}(3)$ \\\hline
      
      \end{tabular}
  \end{center}
\end{table}

$\bullet$ If $B \cong B(3, 2k; 0, \gamma, 0)$ with $k \geq 3$ and $\gamma= 1,2$, then the outer automorphism group of the $\mathcal{F}$-Alperin subgroups can be seen in the table below:
\begin{table}[H]
  \begin{center}
    \caption{Fusion systems over $B(3,2k;0,\gamma,0)$ with $\gamma=1,2$ and $k\geq 3$.}
    \label{tab:table2}
    \begin{tabular}{|c|c|c|c|} 
      \toprule
      $B$ & $E_{0}$ & $V_{0}$ & group realization  \\\hline

      $<\omega>$ & - & $\SL_{2}(3)$ & $\mathcal{F}(3^{2k},2+ \gamma)$ \\\hline
      
      $<\eta\omega>$ & $\SL_{2}(3)$ & - & $3._{\gamma}\PGL_{3}(q)$ \\\hline 
      
      \multirow{2}{*}{$<\eta,\omega>$} & - & $\GL_{2}(3)$ & $\mathcal{F}(3^{2k},2+\gamma).2$\\\cline{2-4}
      
      & $\GL_{2}(3)$ & - &  $3._{\gamma}\PGL_{3}(q).2$ \\\hline
      
      \end{tabular}
  \end{center}
\end{table}
With q a prime power such that $\nu_{3}(q-1)=k-1$.
\smallskip \\

\end{theo}

\section{Application}
We now have enough information on the exotic fusion systems over $B(3,4;0,2,0)$ and $B(3,2k;0,\gamma,0)$, $\gamma =1\text{ or } 2$ in order to check that the conditions in Theorem \ref{Alternative reduction theorem} are satisfied. So supposing that those exotic fusion systems are block realizable then we'll show that they are block realizable by blocks of quasisimple groups. Recall that the reduction theorem of Kessar-Stancu cannot apply in these cases as the considered exotic fusion systems are not reduction simple. Throughout this section we'll use the notations in Lemma~\ref{An F alperin subgroup of B(3,r;0,gamma,0) is}, \\\\


Now we give a lemma that will help us prove Theorem \ref{Condition for ART B(3,4,0,2,0)}.  

\begin{lemm} \label{Aut(v0) in F' for B(3,4,0,2,0) contains SL2(3)}
Let $T$ be a strongly $\mathcal{F}$-closed subgroup of $B=B(3,4;0,2,0)$, let $\tilde {\mathcal{F}}$ be a fusion system containing $\mathcal{F}$, over a $p$-group $\tilde{P}$ containing $P$, and let $\mathcal{F'}$ the fusion system over $T$ such that $\mathcal{F'} \unlhd \mathcal{\tilde{F}}$. Then $T\simeq E_0$ and for all $V \in [V_0]_B, \SL_2(3) \subseteq \Aut_{\mathcal{F'}}(V)$.
\begin{proof}
 For $B=B(3,4,0,2,0)$ and $T$ a normal subgroup of $B$ we have that $T$ intersects non-trivially the center of $B$, and given that $Z(B)=\langle \zeta\rangle$ is of order $3$, we have that $Z(B)$ is contained in $T$. Now $\SL_{2}(3) \leq \Aut_{\mathcal{F}}(V_{0})$ so all the non-trivial elements of $V_{0}$ are $\mathcal{F}$-conjugated. Using that $T$ is $\Ff$-strongly closed, we get that $s$ is in $T$, as it is $\Ff$-cojugated to $\zeta$. Hence using Lamma \ref{normal subgroups of B(3,r,0,gamma,0)} we get that $T=\langle s,\zeta' \rangle=E_{0}$. The subgroups of $E_{0}=\langle s,\zeta' \rangle$ are $V_{0}=\langle s, \zeta \rangle$, $V'_{0}=\langle s \zeta ', \zeta \rangle$, $V''_{0}=\langle s \zeta '^{2}, \zeta \rangle$ and $\langle \zeta, \zeta ' \rangle$. To show that $\Aut_{\mathcal{F'}}(V_{0})$ contains $\SL_{2}(3)$ we want to show that all the subgroups of $V_{0}$ are $\mathcal{F'}$-conjugated. 
 And each of $V_{0}$, $V'_{0}$, $V''_{0}$ have three subgroups of order $3$ other than $Z=\langle \zeta\rangle$, the center of $E_{0}$. For example the subgroups of $V_{0}$ are $C_{0}=\langle  s \rangle$, $C'_{0}=\langle s\zeta \rangle$, $C''_{0}=\langle s\zeta^{2}\rangle$ and $Z$. It's clear that  $C_{0}$, $C'_{0}$ and $C''_{0}$ are $E_{0}$-conjugated, hence $\mathcal{F'}$-conjugated. What remains to show is that these groups are $\mathcal{F'}$-conjugated to $Z$. In fact if we take a map $\alpha$ in $\mathcal{F'}$ that conjugates  $C_{0}$ and  $C'_{0}$, since $\SL_{2}(3) \leq \Aut_{\mathcal{F}}(V_{0}) \leq  \Aut_{\mathcal{\tilde{F}}}(V_{0})$ then there exists a map $\alpha ' $ of order $2$ in $\mathcal{\tilde{F}}$ that conjugates $C_0$ and $Z$ and also conjugates $C'_{0}$ and $C"_{0}$.
 $$\begin{tikzcd}
  C_{0} \arrow{r}{\alpha} \arrow{d}{\alpha'} & C'_{0} \arrow{d}{\alpha'}\\
  Z \arrow{r}{\alpha''} & C''_{0}  
 \end{tikzcd}$$ 
 Now we have  $\mathcal{F'} \unlhd \mathcal{\tilde{F}}$, hence $\alpha''$ is in $\mathcal{F'}$.
 So we have that all subgroups of order $3$ of $V_{0}$ are $\mathcal{F'}$- conjugated and $\Aut_{\mathcal{F'}}(V_{0})$ contains $\SL_{2}(3)$. To conclude, every $V \in [V_0]_B$ is $B$-conjugated, hence $\tilde \Ff$-conjugated to $V_0$. So $\Aut_{\mathcal{F'}}(V)\simeq\Aut_{\mathcal{F'}}(V_0)$.

\end{proof}
\end{lemm}

With the previous Lemma now we can see that that the fusion systems over $B(3,4,0,2,0)$ verify the third condition of Theorem \ref{Alternative reduction theorem} :

\begin{proof}[Proof of Theorem \ref{Condition for ART B(3,4,0,2,0)}]
We need to check that the conditions in the reduction theorem are satisfied for $\mathcal{F}_1=\mathcal{F}(3^{4},3)$ and $\mathcal{F}_2=\mathcal{F}(3^{4},3).2\,$. Conditions 1 and 2 are clearly satisfied. For Condition 3, since we are working in $B(3,4,0,2,0)$ we know that the only proper strongly closed subgroup $T$ in this case is the extra special group $3_{+}^{1+2}$ for which we have a lot of information. In fact, By \cite{RV} we know the outer automorphism group of $T$ in $\mathcal{F'}$ depends on the number of $\mathcal{F'}^{\,cr}$-subgroups of $T$.
By Lemma $\ref{Aut(v0) in F' for B(3,4,0,2,0) contains SL2(3)}$ we showed that $\Aut_{\mathcal{F'}}(V_{0}) \supseteq \SL_{2}(3)$ for every subgroup $V_{0}$ of $T$, and all the subgroups of $T$ are $\mathcal{F'}$-Alperin (since the number of $\Ff'$-Alperin groups could be $2$ or $4$ by \cite{RV} and since we have that all $V \in [V_0]_B$ are $\Ff'$-Alperin we get that $\langle \zeta, \zeta'\rangle$ is also $\Ff'$-Alperin) and $\mid \mathcal{F'}^{\,cr} \mid\geqslant 3$. What we need to see now is whether we have 2 classes or one class of $\Ff'$-Alperin subgroups of $T$. We try to see then if $<\zeta,\zeta'>$ of $T$ is $\mathcal{F'}$ conjugate to the rest of the subgroups of $T$ or not. 
Suppose there is two classes of $\Ff'$-Alperin subgroups and suppose that $V_{0} \overset{\Ff'}{ \sim} V_0'$. If we take  $c_{s_1}$ that conjugates $V_0$ to $V_0'$ and $V_0'$ to $V_0''$ in $\Ff$ i.e. in $\tilde{\Ff}$ we get that $ V_0 \overset{\Ff'}{ \sim} V_0' \overset{\Ff'}{ \sim} V_0''$ since $\Ff' \unlhd \tilde{F}$. Hence they are all $\Ff'$-conjugates; including $\langle \zeta, \zeta'\rangle$; because we can either have two classes $2+2$ or one class of $4$ $\Ff'$-Alperin groups. And by Lemma 4.13 \cite{RV}, we have that $\Out_{\mathcal{F'}}(T)=SD_{16}$ and we get that $\Aut_{\mathcal{F'}}(T)$ is of index $3$ in $\Aut(T)$.\\
In fact, $\Aut_{\mathcal{F'}}(T)= \Inn(T) \rtimes \Out_{\mathcal{F'}}(T)$ and so $\vert\Aut_{\mathcal{F'}}(T)\vert=\vert\Inn(T)\vert\cdot \vert\Out_{\mathcal{F'}}(T)\vert= \vert C_{3}\times C_{3}\vert \cdot \vert SD_{16}\vert= 9\cdot 16$\,.
And, since $\vert\Aut(T)\vert=\vert\Inn(T)\vert \cdot \vert(\Aut(T/\varphi(T))\vert = 9\cdot 3\cdot 2\cdot 8$, we get that $\Aut_{\mathcal{F'}}(T)$ is of index $3$ in $\Aut(T)$.

Given that there is an automorphism of $\Aut_{\tilde {\Ff}}(T)$, of order $3$, outside of $\Aut_{\Ff'}(T)$, given by conj$_{s_1}$, we have that $\Aut_{\mathcal{\tilde{F}}}(T)=\Aut(T)$ and $\vert\Aut_{\mathcal{\tilde{F}}}(T):\Aut_{\mathcal{F'}}(T) \vert =3$. But there is no normal subgroup of index $3$ in $\Aut(T)$, containing $\Inn(T)$. If this would have been true,  using the Frattini map: $\Out_{\mathcal{F'}}(T) \xlongrightarrow{\rho} \GL_{2}(3)$ whose kernel is a $3$-group, we get that $\rho(SD_{16})\in$ Syl$_{2}(\GL_{2}(3))$ is a normal subgroup of $\GL_{2}(3)$ of index $3$, This is  impossible, because there is no normal subgroup of $\GL_{2}(3)$ of index $3$.
 
Hence $\mathcal{F'}\ntrianglelefteq\mathcal{\tilde{F}}$\,. And so there is no fusion system over $T$ normal in $\tilde\Ff$.

\end{proof}

Similarly we give a Lemma that will help with the proof of Theorem \ref{Condition for ART B(3,2k,0,gamma,0)}:

\begin{lemm} \label{Aut(v0) in F' for B(3,2k,0,gamma,0) contains SL2(3)}
Let $T$ be a strongly $\mathcal{F}$-closed subgroup of $B=B(3,2k,0,\gamma,0)$ for $\gamma=1,2$ and $k\geq 3$, let $\tilde {\mathcal{F}}$ be a fusion system containing $\mathcal{F}$, over a $p$-group $\tilde{P}$ containing $P$, and let $\mathcal{F'}$ the fusion system over $T$ such that $\mathcal{F'} \unlhd \mathcal{\tilde{F}}$ . Then $V_0\le T$ and $\Aut_{\mathcal{F'}}(V_{0})$ contains $\SL_{2}(3)$.
\begin{proof}
For $B=B(3,2k,0,\gamma,0)$ with $\gamma=1,2$ and $k\geq 3$, and $T$ strongly $\Ff$-closed we get, similarly to the proof of Lemma $\ref{Aut(v0) in F' for B(3,4,0,2,0) contains SL2(3)}$ that $T$ contains $s$. Hence $T\simeq B(3,2k-1,0,0,0) \simeq \langle s,s_{2}\rangle$, so we have $E_{0}=\langle s_{3}^{3^{k-2}},s_{2}^{3^{k-2}},s \rangle \subset T$. Again, using the same reasoning as in the proof of Lemma \ref{Aut(v0) in F' for B(3,4,0,2,0) contains SL2(3)}, we know that the subgroups $C_{0}=\langle s \rangle,C_{0}'=\langle ss_{2}^{3^{k-2}} \rangle $ and $C_{0}''=\langle ss_{2}^{2.3^{k-2}} \rangle$ of $V_{0}=\langle s_{2}^{3^{k-2}},s\rangle$ are $E_{0}$-conjugated and hence $\mathcal{F'}$-conjugated. We want to show that all the subgroups of order $3$ of $V_{0}$ are $\mathcal{F'}$-conjugated.  What remains to show is that $C_0$, $C_0'$ and $C_0''$ are also $\mathcal{F'}$-conjugated to the center $Z=\langle s_{2}^{3^{k-2}} \rangle$ of $E_0$ (and of $T$).
In fact if we take a map $\alpha$ in $\mathcal{F'}$ that conjugates  $C_{0}$ and  $C'_{0}$, since $\SL_{2}(3) \leq \Aut_{\mathcal{F}}(V_{0}) \leq$  $\Aut_{\mathcal{\tilde{F}}}(V_{0})$ then there exists a map $\alpha ' $ in $\mathcal{\tilde{F}}$ that permutes cyclically $C_0$, $Z$ and $C'_{0}$.
 $$\begin{tikzcd}
  C_{0} \arrow{r}{\alpha} \arrow{d}{\alpha'} & C'_{0} \arrow{d}{\alpha'}\\
  Z \arrow{r}{\alpha''} & C_{0}  
 \end{tikzcd}$$ 
 Now we have  $\mathcal{F'} \unlhd \mathcal{\tilde{F}}$, hence $\alpha''$ is in $\mathcal{F'}$.
 So we have that all subgroups of order $p$ of $V_{0}$ are $\mathcal{F'}$- conjugated and $\Aut_{\mathcal{F'}}(V_{0})$ contains $\SL_{2}(3)$.
\end{proof}
\end{lemm}

We can now show that the fusion systems over $B(3,2k,0,\gamma,0)$, for $\gamma=1,2$, also verify the third condition of Theorem \ref{Alternative reduction theorem} :

\begin{proof}[Proof of Theorem \ref{Condition for ART B(3,2k,0,gamma,0)}]
For $T \unlhd B(3, 2k; 0, \gamma, 0)$, we have that $T\cong B(3, 2k-1; 0, 0, 0)= \langle s,s_{2} \rangle$. Then $\mathcal{F'}$, which is a saturated fusion system over $T$, is one of the fusion systems in Table 6 of \cite{DRV} Lemma 5.9.
Consider the map \\
$\Aut(B(3, 2k-1; 0, \gamma, 0)) \xlongrightarrow{\pi} \Out(B(3, 2k-1; 0, \gamma, 0)) \xlongrightarrow{\rho} \GL_{2}(3)$, where $\pi$ is the canonical projection and $\rho$ is the Frattini map - whose kernel is a $3$-group.

Take the conjugation map $c_{s_{1}} \in \Aut_{\mathcal{F}}(T):$
\begin{center}

$c_{s_{1}}: \langle s,s_{2} \rangle \rightarrow \langle s,s_{2} \rangle$
\\ $s\rightarrow s^{-1}s_{2}$
\\ $s_{2} \rightarrow s_{2}$
\end{center}
Its image via $\rho \circ \pi$ is the lower triangular matrix $\alpha=\left(\begin{smallmatrix}
-1 & 0\\
1 & 1
\end{smallmatrix}\right)$. Now, by Lemma \ref{Aut(v0) in F' for B(3,2k,0,gamma,0) contains SL2(3)}, $\Aut_{\mathcal{F'}}(V_{0})$ contains $\SL_{2}(3)$. 
Hence $\Out_{\mathcal{F'}}(B(3, 2k-1; 0, \gamma, 0))$ contains an element whose image via $\rho$ is the matrix 
$X=\left(\begin{smallmatrix} 
-1 & 0\\
0 & 1
\end{smallmatrix}\right)$. The conjugation of $X$ by $\alpha$ takes us outside $\rho( \Aut_{\mathcal{F'}}(T))$, which is a subgroup of the diagonal matrices in $\GL_{2}(3)$ :\\ 
$^{\alpha}X=\left(\begin{smallmatrix} -1&0\\ 1&1 \end{smallmatrix}\right)\left(\begin{smallmatrix} -1&0\\ 0& 1 \end{smallmatrix}\right) \left(\begin{smallmatrix} 1&0\\-1&1 \end{smallmatrix}\right)=\left(\begin{smallmatrix} 1&0\\1& 1 \end{smallmatrix}\right)$ .
And so there exists an $\alpha \in \rho(\Out_{\mathcal{F}}(T))$ that doesn't normalize $\rho(\Out_{\mathcal{F'}}(T))$ which also means that there exists $\gamma \in \Aut_{\mathcal{F}}(T)\leq \Aut_{\mathcal{\tilde{F}}}(T)$ that doesn't normalize $\Aut_{\mathcal{F'}}(T)$ . 

Hence $\mathcal{F'}\ntrianglelefteq\mathcal{\tilde{F}}$. And so, every fusion system on $T$ is not normal in $\tilde\Ff$.

\end{proof}
 The previous lemma shows that we can apply Theorem \ref{Alternative reduction theorem} to exotic fusion systems over $B(3,4;0,2,0)$ and $B(3,2k;0,\gamma,0)$, for $\gamma=1,2$, $k\ge 2$. Hence if there is a block having the fusion structure of these fusion systems, then there is a block of a quasisimple group with same of structure.

\section*{Acknowledgement}
These results are part of my PhD-Thesis at University Picardie Jules Verne under the supervision of Radu Stancu.\\
My PhD is funded by FEDER and Region Hauts de France.

\bibliographystyle{abbrv}
\bibliography{bibliog.bib}

\end{document}